\documentclass[12pt,twoside]{article}
\title{On
G-$(n,d)$-rings}
\date{}
\usepackage{amsfonts}
\usepackage{amsmath}
\usepackage{amssymb}
\usepackage{latexsym}
\usepackage{graphicx}

\newtheorem{thm}{\bf Theorem}[section]
\newtheorem{cor}[thm]{\bf Corollary}
\newtheorem{lem}[thm]{\bf Lemma}
\newtheorem{prop}[thm]{\bf Proposition}
\newtheorem{defn}[thm]{\bf Definition}
\newtheorem{defns}[thm]{\bf Definitions}
\newtheorem{rem}[thm]{\bf Remark}

\newtheorem{exmp}[thm]{\bf Example}

\catcode`\ç=13
\defç{\c{c}}
\catcode`\é=13
\defé{\'e}
\catcode`\à=13
\defà{\`a}
\catcode`\è=13
\defè{\`e}
\catcode`\â=13
\defâ{\^a}
\catcode`\ù=13
\defù{\`u}
\catcode`\ê=13
\defê{\^e}
\catcode`\î=13
\defî{\^\i}
\catcode`\ô=13
\defô{\^o}

\def\proof{{\parindent0pt {\bf Proof.\ }}}

\def\Ggldim{{\rm G\!-\!gldim}}

\def\Im{{\rm Im}}

\def\Hom{{\rm Hom}}

\def\sup{{\rm sup}}

\newcommand{\cqfd}
{\hspace{1cm}
\rule{2mm}{2mm}%
\medbreak%
\par%
}
\def\1{{\noindent\rm (1)}}
\def\2{{\noindent\rm (2)}}
\def\3{{\noindent\rm (3)}}
\def\4{{\noindent\rm (4)}}
\def\5{{\noindent\rm (5)}}
\begin{document}
\thispagestyle{empty}

\maketitle \vspace*{-1.5cm}
\begin{center}{\large\bf  Najib Mahdou and Khalid Ouarghi}

\bigskip

 \small{Department of Mathematics, Faculty of Science and Technology of Fez,\\ Box 2202, University S. M.
Ben Abdellah Fez, Morocco\\[0.12cm]
mahdou@hotmail.com\\
ouarghi.khalid@hotmail.fr }
\end{center}

\bigskip\bigskip

\noindent{\large\bf Abstract.} The main aim of this paper is to
investigate new class of rings called, for positive integers $n$
and $d$, $G-(n,d)-$rings, over which every $n$-presented module
has a Gorenstein projective dimension at most $d$. Hence we
characterize $n$-coherent $G-(n,0)-$rings. We conclude by various
  examples of $G-(n,d)-$rings.
\bigskip

\bigskip\bigskip


\begin{section}{Introduction}
Throughout this paper all rings are commutative with identity
element and all modules are unital. If $M$ is any $R$-module, we
use $pd_R(M)$, $id_R(M)$ and $fd_R(M)$ to denote,  respectively,
the usual projective, injective and flat dimensions of $M$. It is
convenient to use ``$m$-local" to refer to (not necessarily
Noetherian) ring with a unique maximal ideal $m$.

In 1967-69, Auslander and Bridger \cite{A1,A2} introduced the
  G-dimension for finitely generated modules over Noetherian rings.
Several decades later, this homological dimension was extended, by
Enochs and Jenda \cite{EJ2,EJ1}, to the Gorenstein projective
dimension of modules that are not necessarily finitely generated
and over non-necessarily Noetherian rings. And, dually, they
defined the Gorenstein injective dimension. Then, to complete the
analogy with the classical homological dimension, Enochs,  Jenda
and  Torrecillas \cite{EJT} introduced the Gorenstein flat
dimension.

In the last years, the Gorenstein homological dimensions have
become a vigorously active area of research (see
\cite{BM,LW,EJ2,Rel-hom,HH} for more details). In 2004, Holm
\cite{HH} generalized several results which are already obtained
over Noetherian rings.

The Gorenstein projective, injective and flat dimensions of a
module are defined in terms of resolutions by Gorenstein
projective, injective and flat modules, respectively.

\begin{defns}\cite{HH}\label{def-Gmod}\mbox{}
\begin{enumerate}
\item  An $R$-module $M$ is said to be
Gorenstein projective, if there exists an exact sequence of
projective modules $$\mathbf{P}=\ \cdots\rightarrow P_1\rightarrow
P_0 \rightarrow P^0
         \rightarrow P^1 \rightarrow\cdots$$ such that  $M \cong \Im(P_0
\rightarrow P^0)$ and such that $\Hom_R ( -, Q) $ leaves the
sequence $\mathbf{P}$ exact whenever $Q$ is a projective module.
        \item  The Gorenstein injective modules are defined dually.
        \item An $R$-module $M$ is said to be Gorenstein flat, if there exists an exact sequence of  flat modules
$$\mathbf{F}=\ \cdots\rightarrow F_1\rightarrow F_0 \rightarrow
F^0 \rightarrow F^1 \rightarrow\cdots$$ such that  $M \cong
\Im(F_0 \rightarrow F^0)$ and such that $-\otimes I $ leaves the
sequence $\mathbf{F}$ exact whenever $I$ is an injective module.
\end{enumerate}
\end{defns}

Recently, in  \cite{BM2}, the authors  started the study of
                     the Gorenstein global
                     dimension of a ring $R$, denoted $\Ggldim(R)$, and defined as follows:
$$G-gldim(R)=\sup\{Gpd_R(M)\mid \, M \,\, R-module\}.$$

Let $R$ be a commutative ring, and let $M$ be an $R$-module. For
any positive integer $n$, we say that $M$ is $n$-presented
whenever there is an exact sequence:

$$F_n\longrightarrow F_{n-1}\longrightarrow\ldots \longrightarrow F_0\longrightarrow M \longrightarrow 0$$

 of $R$-modules in which each $F_i$ is a finitely generated free
 $R$-module. In particular, $0$-presented and  $1$-presented $R$-modules are respectively
finitely generated and finitely presented R-modules. We set
$\lambda_R(M)=\sup\{n\mid M$ is $n$-presented$\}$, except that we
set $\lambda_R(M)=-1$ if $M$ is not finitely generated. Note that
$\lambda_R(M)\geq n$ is a way to express the fact that $M$ is
$n$-presented.

Costa (1994) \cite{costa}, introduced a doubly filtered set of
classes of rings in order to categorize the structure of non-
Noetherian rings: for non-negative integers $n$ and $d$, we say
that a ring $R$ is an $(n,d)$-ring if $pd_R(M)\leq d$ for each
$n$-presented $R$-module $M$. $(n,d)$-rings are known rings in
some particular  values of $n$ and $d$. For examples, $R$ is a
Noetherian $(n,d)$-ring,  means that $R$ has global dimension
$\leq d$. $(0,0)$, $(1,0)$, and $(0,1)$-rings are respectively
semi simple, von Neumann regular and hereditary rings (see
\cite[Theorem 1.3]{costa}). According to Costa \cite{costa}, a
ring $R$ is called a  $n$-coherent ring if every $n$-presented
$R$-module is $(n + 1)$-presented. For more results
about $(n,d)$-rings, see for instance, \cite{costa,Mh1,Mh2}.\\

The object of this paper is to extend the ideas of Costa and
introduce a  doubly filtered set of classes of rings called
$G-(n,d)-$rings and  defined as follows:

\begin{defn}\label{dfn2.1}
Let $n$, $d\geq 0$ be  integers.  A ring $R$ is called a
$G-(n,d)$-ring if every $n$-presented $R$-module  has a Gorenstein
projective dimension at most $d$, ($\lambda_R(M)\geq n$ implies
$Gpd_R(M)\leq d$).
\end{defn}

In section 2, we characterize some known rings by
$G-(n,d)-$propriety, for small values of $n$ and $d$. Hence, We
study the transfer of this propriety to some particular
construction rings. Also, we characterize $n-$coherent
$G-(n,0)-$rings. The  section 3 is reserved to give some examples
of $G-(n,d)-$rings for  particular values of $n$ and $d$. Hence,
we give an example of a ring which is a $G-(n,d)-$ring but not
$(n,d)-$ring for any positive integers $n$ and $d$. Also we give
examples of $G-(n,0)-$rings which are not $G-(n-1,d)-$rings, for
$n=2,3$ and for any positive integer $d$.
\end{section}
\begin{section}{Main results}

We start this section by the following Definitions.

\begin{defns}\label{defs2.1}
Let $R$ be a ring.
\begin{enumerate}
    \item $R$ is G-semisimple if  every $R$ module is
    Gorenstein projective (=$R$ is quasi-Frobenius).
    (See \cite[$\S2$]{BMO}).
    \item $R$ is G-Von Neuman regular if  every
    $R$-module is Gorenstein flat (=$R$ is  IF-ring).
    (See \cite[$\S5$]{MT}).
    \item $R$ is G-hereditary if  $G-gldim(R)\leq 1$.
    Also $R$ is a G-Dedekind if it is an integral domain
    G-hereditary. (See \cite[$\S3$]{MT}).
    \item $R$ is G-semi-hereditary if $R$ is coherent
    and every submodule of a flat $R$-module is Gorenstein flat. Also $R$ is a G-pr\"{u}fer if it is an integral domain
    G-semi-hereditary.(See \cite[$\S4$]{MT}).
\end{enumerate}

\end{defns}

 In the next result  we
characterize the  rings  in Definitions \ref{defs2.1} above, with
the $G-(n,d)-$propriety. This Theorem is a generalization of
\cite[Theorem 1.3]{costa}.

\begin{thm}\label{thm2.3}
Let $R$ be a ring. Then:

\begin{enumerate}
    \item $R$ is a $G-(0,0)-$ring if and only if $R$ is
    G-semisimple.
    \item $R$ is a $G-(0,1)-$ring if and only if $R$ is
    G-hereditary.
    \item $R$ is a $G-(0,d)-$ring if and only if $G-gldim(R)\leq d$.
    \item $R$ is a $G-(1,0)-$ring if and only if $R$ is
    G-Von Neuman regular.
    \item If $R$ is  coherent, then $R$ is   $G-(1,1)-$ring if and only if $R$ is
    G-semi-hereditary.
    \item $R$ is a $G-(0,1)-$domain if and only if $R$ is
    G-Dedekind.
    \item If $R$ is  coherent, then $R$ is a   $G-(1,1)-$Domain if and only if $R$ is
    G-pr\"{u}fer.
    \item $R$ is Noetherian, then $R$ is a $G-(n,d)-$ring if and
    only if $G-gldim(R)\leq d$
\end{enumerate}
\end{thm}
\proof $(1)$ Follows from \cite[Proposition 2.1]{BMO}.
$(2),(4),(5),(6)$ and $(7)$ follow respectively from   \cite[
Proposition 3.3, Proposition 5.8, Proposition 4.3, Definition 3.1
and Definition 4.1]{MT}. $(3)$ follows from \cite[Lemma 2.2]{BM2}.
 $(8)$, follows from $(3)$ and since in a Noetherian ring $R$, every finitely generated
$R$-module is infinitely presented.\cqfd

\begin{rem}\label{rem2.2}
1)  An $(n,d)$-ring is a
$G-(n,d)-$ring for  any positive integers $n$ and $d$. The converse is not true in general (see  Example \ref{exmp1}).\\
 2)  $G-(n,d)-$rings are $G-(n',d')-$rings for any $n'\geq n$ and $d\geq
d'$. The converse is not true in general (see Theorem
\ref{thm2.11}).
\end{rem}

Recall that,  for two rings $A \subseteq B$, we say that $A$ is a
module retract of $B$ if there exists an $A$-module homomorphism
$f : B \longrightarrow A$ such that $f/A = id/A$. $f$ is called a
module retraction map. If such  map $f$ exists, $B$ contains $A$
as an $A$-module direct summand.
\begin{prop}\label{thmretract}
Let  $A$ be a  subring retract of $R$, ($R=A\oplus_A E$), such
that $E$ is a flat $A$-module and $G-gldim(A)$ is finite. If $R$
is a $G-(n,d)-$ring, then $A$ is a $G-(n,d)-$ring too.
\end{prop}
\proof

Let $M$ be an $A$-module $n$-presented. Since $R$ is a flat
$A$-module, $M\otimes_A R$ is an $R$-module $n$-presented, and by
hypothesis $Gpd_R(M\otimes_A R)\leq d$. Then,  $Gpd_A(M)\leq d$
from \cite[Proposition2.4]{MO} .\cqfd

Let  $A$ be a ring and let $E$ be  an $A$-module. The trivial ring
extension of  $A$ by $E$ is the ring $R := A \propto E$ whose
underlying group is $A \times E$ with multiplication given by
$(a,e)(a',e') = (aa',ae'+a'e)$. These extensions have been useful
for solving many open problems and conjectures in both commutative
and non-commutative ring theory. See for
instance,\cite{Glaz,Hu,Kabbajmahdou2,Mh1,Mh2}.\\

 A direct application of
Proposition \ref{thmretract} is the following Corollary.

\begin{cor}\label{cor Trivial}
Let $A$ be a ring and let $E$ be an $A$-module such that $E$ is
flat
 and $G-gldim(A)$ is finite. If $R=A\propto E$ is a
$G-(n,d)-$ring, then $A$ is a $G-(n,d)-$ring too.
\end{cor}

In the next result we study the transfer of the
$G-(n,d)-$propriety to the polynomial ring.

\begin{thm}\label{R[X]}
Let $R$ be a ring and let $X$ be an indeterminate over $R$.

\begin{enumerate}
    \item Suppose that $G-gldim(R)$ is finite. If $R[X]$ is a
$G-(n,d)-$ring, then $R$ is a $G-(n,d)-$ring too.
    \item If $R$ is a $G-(n,d)-$ring which is not
    $G-(n,d-1)$-ring, then $R[X]$ is not a $G-(n,d)-$ring.
    \item Suppose that $G-gldim(R)$ is finite. If $R[X]$ is a
    $G-(n,d)-$ring, then $R$ is a $G-(n,d-1)$-ring.
\end{enumerate}
\end{thm}

\proof \begin{enumerate}
    \item  Let $M$ be an $R$-module such that
$\lambda_R(M)\geq n$.
 Since $R[X]$ is a free $R$-module, we have $\lambda_{R[X]}(M[X])\geq
n$, and by hypothesis $Gpd_{R[X]}(M[X])\leq d$. From \cite[Lemma
2.8]{BMpoly} $Gpd_R(M)\leq d$, and $R$ is a $G-(n,d)-$ring as
desired.
    \item Since $R$ is a $G-(n,d)$-ring which is  not $G-(n,d-1)$-ring,
    there exists an $R$-module $M$ such that $\lambda_{R}(M)\geq
n$ and $Gpd_R(M)=d$. Then, from \cite[Theorem]{HH}, it is easy to
see   that there exists a free $R$-module $F$ such that
$Ext_R^{d}(M,F)\neq 0$. on the other hand, $M$ is also an
$R[X]$-module via the canonical morphism: $R[X]\longrightarrow R$.
Hence, from \cite[Lemma 9.29]{Rot}, there exists an exact sequence
of $R[X]$-modules: $$0\longrightarrow M[X]\longrightarrow
M[X]\longrightarrow M\longrightarrow 0,$$ from which we conclude
that $\lambda_{R[X]}(M)\geq \lambda_{R[X]}(M[X])$. But since
$R[X]$ is a flat $R$-module we see that $\lambda_{R[X]}(M[X])\geq
\lambda_{R}(M)\geq n$, and we have $\lambda_{R[X]}(M)\geq n$.
Then, \cite[Theorem 9.37]{Rot} shows that:
$$Ext_{R[X]}^{d+1}(M, F[X])\cong Ext_{R}^{d+1}(M, F)\neq 0.$$
It follows from \cite[Theorem 2.20]{HH} that $Gpd_{R[X]}(M)\geq
d$. Finally,  $R[X]$ is not a $G-(n,d)-$ring as desired.
    \item . Follows from (1) and (2) of the same Theorem.\cqfd
\end{enumerate}

In the next Theorem we study the transfer of the
$G-(n,d)-$propriety to finite direct product of rings.

\begin{thm}\label{thm2.5}
Let  $R=R_1\times R_2...\times R_m$  be a finite direct product of
rings. If  $R$ is $G-(n,d)-$ring, then
   $R_i$ is a $G-(n,d)-$ring for each $i=1,...,m$.\\
The converse is true if $sup\{G-gldim(R_i)\mid i=1,...,m\}$ is
finite.
\end{thm}
To prove this Theorem we need the following Lemma.

\begin{lem}\label{lemn-presented}
Let  $R=R_1\times R_2...\times R_m$  be a finite direct product of
rings and let $n\geq 0$ be an integer. Then,  $M=\oplus_{i}M_i$ is
$n-$presented $R$-module if and only if, $M_i$ is $n$-presented
$R_i-$module for each $i=1,...,m$.
\end{lem}
 \proof
Follows from \cite[Corollary 2.6.9]{BK}.\cqfd
\noindent\textbf{Proof of Theorem \ref{thm2.5}.}
  Let $M_i$ be an $R_i$-module such that $\lambda_{R_i}(M_i)\geq n$,
then, from Lemma \ref{lemn-presented} above, we have
$\lambda_{R}(\oplus_{i}M_i)\geq n$ and by hypothesis
$Gpd_R(\oplus_{i}M_i)\leq d$. Hence, from \cite[Lemma 3.2]{BMpoly}
$Gpd_{R_i}(M_i)\leq d$.\\
Conversely, suppose that  $sup\{G-gldim(R_i)\mid i=1,...,m\}$ is
finite and let $M=M_1\oplus...\oplus M_m$ be an $R$-module
$n$-presented. Then for each $i$, $M_i$ is $n$-presented
$R_i$-module by Lemma \ref{lemn-presented}.And  from the
hypothesis we have $Gpd_{R_i}(M_i)\leq d$. Hence, from \cite[Lemma
3.3]{BMpoly}, $Gpd_R(M)\leq \sup\{Gpd_{R_i}(M_i)\mid
i=1,...,m\}\leq d$.\cqfd

The next result shows that a $G-(n,d)-$ring has grade at most $d$.
This Theorem is a generalization of \cite[Theorem 1.4]{costa}.

\begin{thm}\label{thm2.7}
Let $R$ be a $G-(n,d)-$ring. Then $R$ contains no regular sequence
of length $d+1$.
\end{thm}
\proof Let $x_1,...,x_t$  be a regular sequence in $R$, where
$I=\sum_{i=1}^{t}Rx_i\neq R$. Then, the Koszul complex defined by
$\{x_1,...,x_t\}$ is a finite free resolution of $R/I$ and hence
$R/I$ is $n$-presented for every $n$. Then, since $R$ is a
$G-(n,d)-$ring, we have $Gpd_R(R/I)\leq d$. But it's well known
from \cite[Exercice 1, page 127]{Kap}, that $Gpd_R(R/I)=t$ . Hence
$t\leq d$.\cqfd

In the next result we study "when the $G-(n,d)-$propriety is
local" of the $G-(n,d)-$propriety.
\begin{prop}\label{prop2.9}
Let $R$ be a ring with $G-gldim(R)$ is finite and let $n$ and $d$
be positive integers such that, $d\leq n-1$. If $R$ is  locally a
$G-(n,d)-$ring, then $R$ is also a $G-(n,d)-$ring.
\end{prop}
 To prove this Theorem we need this Lemma.
\begin{lem}\cite[Lemma 3.1]{costa}\label{lemmacosta}
Let $M$ be an $R$-module, and let $S$ be a multiplicative subset
of  system in $R$. If $M$ has a finite $n$-presentation, then:
$$S^{-1}Ext_R^i(M,N)\cong Ext_{S^{-1}R}^i(S^{-1}M,S^{-1}N)$$
for all $0\leq i\leq n-1$, and $S^{-1}Ext_R^n(M,N)$ is isomorphic
to some  submodule of $Ext_{S^{-1}R}^n(S^{-1}M,S^{-1}N)$.
\end{lem}
\noindent\textbf{Proof of Proposition \ref{prop2.9}.} Let $M$ be
an $n$-presented $R$-module  and let $m$ be a maximal ideal of
$R$. Then $M_m$ is an  $n$-presented $R_m$-module. Let $P$ be a
projective $R$-module, then
$(Ext_R^i(M,P))_m=Ext_{R_m}(M_m,P_m)=0$, and from
\cite[Theorem3.80]{Rot}, $Ext_R^i(M,P)=0$ for all $0\leq i\leq n$.
Hence, $Gpd_R(M)\leq d$.\cqfd

Now we give our main result of this section  in which we give a
 characterization of $n$-coherent $G-(n,0)-$ring.

\begin{thm}\label{thm2.10}
Let $R$ be an $n$-coherent ring. Then the following conditions are
equivalent.
\begin{description}
    \item[A)]$R$ is a $G-(n,0)-$ring.
    \item[B)] The following conditions hold:

\begin{enumerate}
    \item Every finitely generated ideal of $R$ has a nonzero
    annihilator.
    \item For each infinitely presented $R$-module $M$
         $Gpd_R(M)<\infty$.
    \item For every finitely generated Gorenstein projective
        submodule $G$ of a finitely generated projective $R$-module $P$, $P/G$ is
        Gorenstein projective.
\end{enumerate}
\end{description}
\end{thm}
\proof

To prove this Theorem we need the following Lemma.
\begin{lem}[\cite{bass}, Theorem 5.4]\label{lemmabass}
The following assertions are equivalent for a ring $R$:
\begin{enumerate}

    \item Every finitely generated projective submodule of a
    projective $R$-module $P$ is a direct summand of $P$.
    \item Every finitely generated proper ideal of $R$ has a nonzero
    annihilator.
\end{enumerate}
\end{lem}
\noindent\textbf{Proof of Theorem \ref{thm2.10}.}
$(A)\Longrightarrow (B).$ The condition (2) is obvious. Now we
prove  $(1).$ Let $P$ be a finitely generated submodule of $Q$ and
    both $P$ and $Q$ are projective. Let $Q'$ be a projective $R$-module such that $Q\oplus_R Q'=F_0$ is a free $R$-module. Then there exists an  exact sequence:
    $$0\longrightarrow P\longrightarrow F_0\longrightarrow Q/P\oplus_R Q'\longrightarrow 0 \qquad\qquad (\ast)$$

      On the other hand, since $P$ is finitely
    generated projective $R$-module, there exists a finitely
    generated free submodule $F_1$ of $F_0$   such that $P\subseteq F_1$ and $F_0=F_1\oplus_R
    F_2$. Thus,   we see easily that $P$ is
    infinitely presented and from the exact sequence: $$0\longrightarrow P\longrightarrow F_1\longrightarrow F_1/P\longrightarrow
    0$$
    $pd_R(F_1/P)\leq 1$ and $F_1/P$ is also infinitely presented, and by hypothesis  $F_1/P$ is Gorenstein projective, then it is
    projective. Consider the following pushout diagram:

$$
\begin{array}{ccccccccc}
   &  &  &  & 0 &  & 0 &  &  \\
   &  &  &  & \downarrow &  & \downarrow &  &  \\
  0 & \rightarrow & P & \rightarrow & F_1 & \rightarrow & F_1/P & \rightarrow & 0 \\
   &  & \shortparallel &  & \downarrow &  & \downarrow &  &  \\
  0 & \rightarrow & P & \rightarrow & F_0 & \rightarrow & Q/P\oplus_R Q' & \rightarrow & 0 \\
   &  &  &  & \downarrow &  & \downarrow &  &  \\
   &  &  &  & F_2 & = & F_2 &  &  \\
   &  &  &  & \downarrow &  & \downarrow &  &  \\
   &  &  &  & 0 &  & 0 &  &  \\
\end{array}%
$$

Since $F_2$ and $F_1/P$ are projective, the exact sequence
$(\ast)$ splits and $Q\oplus_R Q'\cong P\oplus_R Q/P\oplus_R Q'$.
Then
$Q\cong P\oplus_R Q/P$ as desired.\\
To finish the proof of the first implication, it remains to prove
that the condition $(3)$ holds. Let $G$ be a finitely generated
submodule of a finitely generated projective $R$-module $P$.
Consider the exact sequence:$$0\longrightarrow G \longrightarrow
P\longrightarrow P/G\longrightarrow 0.$$ It follows that
$Gpd_R(P/G)\leq 1$ and from \cite[Theorem 2.10]{HH}, there exists
an exact sequence of $R$-modules : $$0\longrightarrow
K\longrightarrow H\longrightarrow P/G\longrightarrow 0\quad
(\star)$$ where $K$ is projective and $H$ is Gorenstein
projective. Combine $(2)$ of this Theorem with \cite[Proposition
5.17]{MT}, we conclude that $K$ is a direct summand of $H$ and the
exact sequence $(\star)$ splits. Then $P/G$ is
Gorenstein projective as  direct summand of $H$.\\

$(B)\Longrightarrow (A)$

  Let $M$ be an $n$-presented $R$-module. Since $R$ is
    $n$-coherent, $M$ is infinitely presented and $Gpd_R(M)$ is
    finite. Let $Gpd_R(M)=d$, then we have the exact sequence of $R$-modules:
    $$0\longrightarrow G\stackrel{u_{d-1}}{\longrightarrow} P_{d-1}\stackrel{u_{d-2}}{\longrightarrow} P_{d-2}\cdots \longrightarrow
    P_1
    \stackrel{u_{1}}{\longrightarrow} P_{0}\stackrel{u_{0}}{\longrightarrow} M\longrightarrow 0 $$
    where $P_i$ is finitely generated projective for each $i$ and $G$ is
    Gorenstein projective. Then we have the exact sequences of $R$-modules:
    $$0\longrightarrow G(=\ker(u_{d-1}))\longrightarrow P _{d-1} \longrightarrow Im(u_{d-1})\longrightarrow 0$$
    $$0\longrightarrow Im(u_i)(=\ker(u_{i-1}))\longrightarrow P _{i-1} \longrightarrow Im(u_{i-1})\longrightarrow 0\quad \hbox{for } i=2,...,d-1$$
    $$0\longrightarrow Im(u_1)(=\ker(u_{0}))\longrightarrow P _{0} \longrightarrow Im(u_{0})=M\longrightarrow 0$$
    Then, by hypothesis and  since $G$ is a finitely generated Gorenstein projective submodule of a
    projective $R$-module $P_{d-1}$, we have  $Im(u_{d-1})\cong
    P_{d-1}/G$ is a  finitely generated Gorenstein projective $R$-module.
    Thus, by
    induction, we conclude that $M=Im(u_0)$ is a finitely
    generated Gorenstein projective $R$-module and this completes
    the proof.\cqfd

In the next Proposition we study the relation between
$G-(n,d)-$rings and $G-(n,0)-$rings.
\begin{prop}\label{prop2.8}
Let $R$ be a $G-(n,d)-$ring. Then $R$ is a $G-(n,0)-$ring if and
only if  $Ext_R(M,K)=0$ for every $n$-presented $R$-module $M$ and
every $R$-module $K$ with $pd_R(K)=Gpd_R(M)-1$.
\end{prop}
\proof

$\Longrightarrow)$ Obvious.\\
$\Longleftarrow)$ Let $M$ be an $n$-presented $R$-module. Since
$R$ is a $G-(n,d)-$ring, we have $Gpd_R(M)\leq d$. And from
\cite[Theorem 2.10]{HH} there exists an exact sequence:
$$0\longrightarrow K\longrightarrow G\longrightarrow
M\longrightarrow 0\qquad\qquad(\star)$$ where $G$ is Gorenstein
projective and $pd_R(K)=Gpd_R(M)-1$. By hypothesis $Ext_R(M,K)=0$
and the exact sequence $(\star)$ splits. Then by \cite[Theorem
2.5]{HH} $M$ is Gorenstein projective as  a direct summand of
$G$.\cqfd

\end{section}

\begin{section}{Examples}

In this section, we construct a class of $G-(2,0)-$rings
(respectively, $G-(3,0)-$rings) which are not $(1,d)-$rings
(respectively, not $G-(2,d)-$rings) for every integer $d\geq 1$.
Also we give an example of ring which is $G-(n,d)-$ring and not
$(n,d)-$ring for every integers $n,d\geq 0 $.

In the next result we give an example of ring which is
$G_(2,0)-$ring but not $G-(1,d)-$ring. Also we give an example of
ring which is $G-(2,d)-$ring and neither $G-(2,d-1)-$ring nor
$G-(1,d)-$ring for any integer $d\geq 0$. This Theorem is a
generalization  \cite[Theorem 3.4]{Mh1}.

\begin{thm}\label{thm2.11}
 Let $K$ be a field and let $E(\cong K^{\infty})$  be a $K$-vector space with infinite
 rank. Let $R:=K\propto E$ the trivial ring extension of $K$ by
 $E$. Then:
 \begin{enumerate}
    \item $R$ is a $G-(2,0)-$ring.
    \item  $R$ is not $G-(1,d)-$ring, for every positive integer
    $d$.
    \item  Let $S$ be a Noetherian ring with $Ggldim(S)=d$. Then,
    $T=R\times S$ is a $G-(2,d)-$ring but neither $G-(1,d)-$ring nor
    $G-(2,d-1)-$ring
 \end{enumerate}

\end{thm}
\proof

\begin{enumerate}
    \item $R$ is a $G-(2,0)-$ring since it is a $(2,0)-$ring from
    \cite[Theorem  3.4]{Mh1}.
    \item Let $d$ be a positive integer, we have to  prove that $R$
    is not a $G-(1,d)-$ring.  $M=0\propto E$ is the maximal
    ideal of $R$ and let $(0,e_i)_{i\in}$ be a set of generators of $M$. Consider the exact sequence of $R$-modules:
    $$0\longrightarrow M^{(I)}\longrightarrow R^{(I)}\longrightarrow M\longrightarrow 0$$
    from this exact sequence we conclude that $Gpd_R(M)=0$ or $Gpd_R(M)=\infty$.
    Suppose that $Gpd_R(M)=0$ and let $J=R(0,f)\cong 0\propto K$ be a principal
    ideal of $R$. $J$ is a direct summand  of $M$, then  $Gpd_R(J)=0$. Consider the exact sequence of $R$-modules:
    $$0\longrightarrow ker(u)\longrightarrow R \stackrel{u}{\longrightarrow} J\longrightarrow 0$$
    where
    $u((a,e))=(a,e)(0,f)=(0,af)$.
    Then, $ker(u)=\{(a,e)\in R\mid
    af=0\}$. Easily we can conclude that
    $ker(u)=M$. Then, $M\cong R/J \cong K$, hence  $K$ is a Gorenstein projective $R$-module.
    In particular $Ext_R(K,R)=0$, and  $R$ is self-injective($0=id_K(E)=id_R(R)$) from \cite[Proposition
    4.35]{FG}, contradiction. Indeed, $R$ is not self-injective since $Ann_R(ann_R((J))=M\neq
    J$ and from \cite[Corollary 1.38]{NY}. Then $Gpd_R(M)=\infty$. On the other hand, $R/J$ is
    $1$-presented $R$-module and $Gpd_R(M)=Gpd_R(R/J)=\infty$.
    Finally,
    $R$ is not a $G-(1,d)-$ring for each positive integer $d$.
    \item Follows from Theorem \ref{thm2.5} and Theorem
    \ref{thm2.3} (8). \cqfd
\end{enumerate}

Next we give an example of  a $G-(n,d)-$ring which is not an
$(n,d)$-ring, for every positive integers $n$ and $d$.

\begin{exmp}\label{exmp1}
Let $K$ be a field and $R=K\propto K$ the trivial ring extension
of $K$ by $K$. Then $R$ is a $G-(n,d)-$ring but not an
$(n,d)$-ring, for every positive integers $n$ and $d$.
\end{exmp}
\proof  From \cite[Theorem 3.7]{BMO}, $R$ is a $G-(0,0)-$ring
(=quasi-Frobenius), then, from Remark \ref{rem2.2} $R$ is a
$G-(n,d)-$ring for every positive integers $n$ and $d$. And it
follows from \cite[Example 3.4]{Mh2} that $R$ is not a
$(n,d)$-ring.\cqfd

 Next result generates  an example of ring which is
$G_(3,0)-$ring but not $G-(2,d)-$ring for every integer $d\geq 0$.
Also we give an example of  $G-(3,d)-$ring which is neither
$G-(3,d-1)-$ring nor $G-(2,d)-$ring.

\begin{thm}\label{thm(3,0)}
Let $(A,M)$ be a local ring and let $R=A\propto A/M$ be the
trivial ring extension of $A$ by $A/M$. Then:
\begin{enumerate}
    \item If $M$ is not finitely generated, then $R$ is a
    $G-(3,0)-$ring.
    \item If $M$ contains a regular element, then $R$ is not a
    $G-(2,d)-$ring, for every integer $d\leq 0$.
    \item Let $S$ be a Noetherian ring with $G-gldim(S)=d$ for some  integer $d\geq0$. Then,
    $T=R\times S$  is a $G-(3,d)-$ring which is neither $G-(2,d)$ nor
    $G-(3,d-1)-$ring.
\end{enumerate}
\end{thm}
\proof

\begin{enumerate}
    \item Follows from \cite[Theorem 1.1]{Kabbajmahdou}.
    \item Suppose that $M$ contains a regular element. Consider the exact
    sequence of $R$-modules:
    $$0\longrightarrow M\propto A/M\longrightarrow R\longrightarrow R/(M\propto A/M)\longrightarrow 0.\qquad (\star)$$
    We claim that $Gpd_R(R/(M\propto
    A/M))=\infty$. Deny, $Gpd_R(R/(M\propto
    A/M))$ is finite. From the exact sequence $(\star)$ and \cite[Proposition  2.18]{HH} we have:
    $$Gpd_R(M\propto A/M)+1= Gpd_R(R/(M\propto
    A/M))\qquad (1)$$

Let $(x_i)_{i\in I}$ be a set of generators of $M$ and $R^{(I)}$
be a free $R$ module. Consider the exact sequence   of
$R$-modules:
$$0\longrightarrow ker(u) \longrightarrow R^{(I)}\oplus_R R \stackrel{u}{\longrightarrow} M\propto A/M\longrightarrow 0$$
where $u((a_i,e_i)_{i\in I},(b_0,f_0))=\sum_{i\in
I}(a_i,e_i)(x_i,0)+(b_0,f_0)(0,1)=(\sum_{i\in I}(a_ix_i,b_0)$,
since $x_i\in M$ for each $i\in I$. Hence
$$ker(u)=(U\propto (A/M)^{(I)})\oplus_R (M\propto A/M)$$ where
$U=\{(a_i)_{i\in I}\in A^{(I)}\mid \sum_{i\in I}a_ix_i=0\}$.
Therefore, we have the isomorphism of $R$-modules:
$$M\propto A/M\cong [R^{(I)}/(U\propto (A/M)^{(I)})]\oplus_R [R/(M\propto A/M)].$$
Hence from \cite[Proposition 2.19]{HH}, we have:
$$Gpd(R/(M\propto A/M))\leq Gpd_R(M\propto A/M)\qquad (2).$$
It follows from $(1)$ and $(2)$ that $Gpd(R/(M\propto A/M))=
Gpd_R(M\propto A/M)=\nolinebreak\infty.$\\
Now from the exact sequence of $R$-modules:
$$0\longrightarrow M\propto A/M \longrightarrow R \longrightarrow 0\propto A/M\longrightarrow 0,$$
we conclude that  $Gpd_R(0\propto A/M)=\infty$. On the other hand,
let $m\in M$ be a regular element and $J=R(m,0)$ an ideal of $R$.
Consider the exact sequence of $R$ modules:
$$0\longrightarrow ker(v)\longrightarrow R\stackrel{v}{\longrightarrow} J\longrightarrow 0,$$
where $v((b,f))=(b,f)(m,0)=(bm,mf)$. Since $m$ is a regular
element, we have $ker(v)=0\propto A/M$. Therefore, it follows that
$Gpd_R(J)=Gpd_R(0\propto A/M)=\infty$. On the other had, $0\propto
A/M$ is a finitely generated ideal of $R$, hence $J$ is a finitely
presented ideal of $R$. Finally, the exact sequence of
$R$-modules:
$$0\longrightarrow J\longrightarrow R\longrightarrow R/J\longrightarrow
0,$$ shows that $\lambda_R(R/J)\geq 2$ and $Gpd_R(R/J)=\infty$.
Then $R$ is not a $(2,d)-$ring for each positive integer $d$.
    \item Follows from Theorem \ref{thm2.5} and Theorem
    \ref{thm2.3} (8). \cqfd
\end{enumerate}

\end{section}

\end{document}